\documentclass{article}
\usepackage{graphicx} 
\usepackage{amsmath, amssymb, amsfonts}
\usepackage{bm}
\usepackage{geometry}
\usepackage[english]{babel}
\usepackage[sort&compress]{natbib}
\newtheorem{theorem}{Theorem}[section]
\newtheorem{definition}{Definition}[section]
\newtheorem{lemma}{Lemma}[section]

\newtheorem{remark}{Remark}

\allowdisplaybreaks[4]

\title{Limit Law for the Maximum Interpoint Distance of High Dimensional Dependent Variables}
\author{Guowei Yan and Long Feng\\
School of Statistics and Data Science, KLMDASR, LEBPS, and LPMC, Nankai University}
\date{}

\begin{document}

\maketitle

\begin{abstract}
In this paper, we considier the limiting distribution of the maximum interpoint Euclidean distance $M_n=\max _{1 \leq i<j \leq n}\left\|\boldsymbol{X}_i-\boldsymbol{X}_j\right\|$, where $\boldsymbol{X}_1, \boldsymbol{X}_2, \ldots, \boldsymbol{X}_n$ be a random sample coming from a $p$-dimensional population with dependent sub-gaussian components. When the dimension tends to infinity with the sample size, we proves that $M_n^2$ under a suitable normalization asymptotically obeys a Gumbel type distribution. The proofs mainly depend on the Stein-Chen Poisson approximation method and high dimensional Gaussian approximation.
\end{abstract}


\section{Introduction}

Let $\bm{X}_1 , \dots , \bm{X}_n$ be identically distributed copies of $p$-dimensional random vector $\bm{X}$ with mean $\bm{\mu}$ and covariance matrix $\mathbf{R}$. Write $\bm{X}_i=\left(x_{i,1},x_{i,2},\dots,x_{i,p} \right)^T$ for each $i=1,2,\dots,n$, and $\bm{X}^{\left(j\right)}=\left(x_{1,j},x_{2,j},\dots,x_{n,j} \right)^T$ for each $j=1,2,\dots,p$. In this paper, we are concerned the asymptotic distribution of the maximum interpoint distance \\
\begin{equation}
M_n=\max\left \{ \left | \bm{X}_i-\bm{X}_j  \right |;1\le i<j\le n \right \} ,
\end{equation}
where $\left|\cdot\right|$ is the Euclidean norm on $\mathbb{R}^p$. It’s clear that the greatest distance is achieved when two vectors, $\bm{X}_i$ and $\bm{X}_j$, which are nearly the longest in length, are almost diametrically opposed. Furthermore, the statistic $M_n$ has potential uses in outlier detection, as discussed by \citet{ref1}.

In recent literature, it has been demonstrated by various researchers that the asymptotic distribution of $M_n$ primarily falls into one of three distinct extreme-value distributions, namely the Gumbel, Weibull, and Fréchet distributions.
For the univariate case, the asymptotic distribution of $M_n$ has been well established, see \citet{ref2} for more information. Specially, \citet{ref3} show that the asymptotic distribution of $M_n$ are sum of two independent Gumble distributions under the assumption that $\bm{X_1},\cdots,\bm{X_n}$ are i.i.d symmetric real-valued random variables.
For the multivariate case, many literature investigated the limiting behavior of $M_n$ under different distribution assumptions. For the standard multivariate normal distribution, \citet{ref1} showed that
\begin{align*}
P\left(\sqrt{2 \log n}\left[M_n-2 \sqrt{2 \log n}-\frac{(p-3) / 2 \log \log n+\log \log \log n+a}{\sqrt{2 \log n}}\right] \leq x\right)\to e^{-e^{-x}},
\end{align*}
where $a=\log \frac{(p-1) 2^{(p-7) / 2}}{\Gamma(p / 2) \sqrt{\pi}}$. \citet{ref4},\citet{ref3} extended the above result to the symmetric Kotz population, spherically symmetric type distribution, respectively. A multitude of findings regarding the maximum interpoint Euclidean distance have been derived when the population $\bm{X}$ is associated with different kinds of multivariate distributions. These findings can be broadly classified into two groups based on whether the distribution of $\bm{X}$ has a bounded or unbounded support. For the unbounded support case, \citet{ref5} showed that the limiting distribution of $M_n$ is not of the Gumble type if the distribution of $\bm{X}$ is power-tailed spherically decomposable. \citet{ref6} consider the elliptical distributions with unbounded support. In the case of bounded support, \citet{ref7} provide a comprehensive review of the current body of literature.

All the aforementioned results have been determined for a population with a fixed dimension. Therefore, it would be intriguing to explore a similar geometric phenomenon in a high-dimensional context. When the population $\bm{X}$ have independent components and each component obeys the sub-exponential distribution, \citet{ref7} showed that the limiting distribution of $M_n$ is
\begin{align}\label{expd}
4 \sqrt{\log n}\left(\frac{M_n^2-2 p}{\sqrt{2\left(\kappa_4+1\right) p}}-2 \sqrt{\log n}+\frac{\log \log n}{4 \sqrt{\log n}}\right) \stackrel{d}{\rightarrow} \xi,
\end{align}
where $\xi$ is a random variable whose distribution is $F_{\xi}(x)=e^{-K e^{-x / 2}}$ with $K=\frac{1}{4 \sqrt{2 \pi}}$ and $\kappa_4=E(X_{ij}^4)$. In practice, the independent assumption is too restrictive. So we will consider the limiting distribution of $M_n$ under a special dependent structure in high dimensional settings.

In this paper, we will consider $\bm X=\bm \mu+\mathbf{R}^{1/2}\bm \epsilon$ where $\bm \epsilon$ are i.i.d sub-gaussian variables and $\mathbf{R}$ is a Toeplitz matrix with a sequence $\left \{ r_k \right \} ^\infty _{k=1}$, that is,
\begin{equation}
\begin{gathered}
\mathbf{R}=  \left ( r_{\left | i-j \right | } \right ) _{p\times p}=
\begin{pmatrix}
 1 & r_1 & \dots & r_{p-1}\\
 r_1 & 1 & \dots & r_{p-2}\\
 \vdots  & \vdots  & \; & \vdots \\
 r_{p-1} & r_{p-2} & \dots & 1
\end{pmatrix} ,
\end{gathered}
\end{equation}
where $r_k>0$ for any $k=1,2,\dots$. According to the Chen–Stein Poisson approximation method and high dimensional Gaussian approximation, we also show that the limiting distribution of $M_n$ still is Gumble type.

\section{Main results}

\begin{definition}(mixing coefficient)
    Let $\left\{\bm{X}_t\right\}$ be a random sequence. Denote by $\mathcal{F}_{-\infty}^u $ and $\mathcal{F}^{\infty}_u$ the $\sigma$-fields generated respectively by $\left\{\bm{X}_t\right\}_{t\le u}$ and $\left\{\bm{X}_t\right\}_{t\ge u}$. The $\alpha$-mixing coefficient and $\rho$-mixing coefficient at lag $k$ of the sequence $\left\{\bm{X}_t\right\}$ is defined as
    \begin{align*}
        \alpha_n\left(k\right)&:=\sup_{t}\sup_{A\in \mathcal{F}_{-\infty}^t, B\in \mathcal{F}^{\infty}_{t+k} } \left | \mathbb{P}\left ( AB \right )- \mathbb{P}\left ( A \right ) \mathbb{P}\left ( B \right )  \right | \\
        \rho_n\left(k\right)&:=\sup_{t}\sup_{f\in L^2\left(\mathcal{F}_{-\infty}^t\right), g\in L^2\left(\mathcal{F}^{\infty}_{t+k}\right) }\left|\mathrm{corr}\left(f,g\right)\right|
    \end{align*}
    We say the sequence $\left\{\bm{X}_t\right\}$ is $\alpha$-mixing or $\rho$-mixing if $\alpha_n\left(k\right)\rightarrow 0$ or $\rho_n\left(k\right)\rightarrow 0$ as $k\rightarrow\infty$. Clearly $\alpha_n\left(k\right)\le \rho_n\left(k\right)$ for any $k$.
\end{definition}

Next, we present the conditions required by our theorem.

\begin{itemize}
\item[(C1)]  Suppose the stationary sequence $X_{i1},\cdots,X_{ip}$ for any $i=1,\dots,n$ satisfies the strong mixing condition and the corresponding $\rho$-mixing coefficient $\rho_n\left(k\right)$ satisfies $\rho_n \left ( k \right ) \le K_1\exp \left ( -K_2k^{\gamma } \right ) $ for any $k\ge 1$ for some universal constants $K_1>1$ , $K_2>0$ and $\gamma>0$.
\item[(C2)] Let $ \nu =\frac{2+4\gamma }{3\gamma }\vee \frac{7}{3} $ we assume $\left ( \log n \right ) ^\nu=o\left ( p^{1/9} \right ) $ where $\gamma$ has been defined in Condition (C1).
\end{itemize}

\begin{remark}
According to \citet{ref8}, we know that, for integer time
\begin{equation*}
    \rho\left(k\right)=\inf_{\phi}\operatorname*{ess\,sup}_{\lambda}\left [ \left | f\left ( \lambda  \right )-e^{i\lambda k}\phi\left ( e^{-i\lambda } \right )  \right | \frac{1}{f\left ( \lambda  \right ) } \right ],
\end{equation*}
where $\inf_\phi$ is taken over all functions $\phi\left(z\right)$ that can by analytically continued inside the unit disc and the spectral density function $f\left(\lambda\right)$ is
\begin{equation}
    f\left ( \lambda  \right ) =\frac{1}{\pi}\sum_{k=0}^{\infty }r_ke^{-ik\lambda },
\end{equation}
where $r_0=1$.  So we know that $\rho_n$ can be bounded by $\left \{ r_k \right \} ^\infty _{k=1}$.
\end{remark}

\begin{remark}
    There are two simple examples of $\mathbf{R}$ satisfying Condition (C1). \\

    $\mathbf{R}$ is the covariance matrix of the  auto-regressive model $AR\left(1\right)$, that is,\\
    \begin{equation*}
        \mathbf{R}=\left(r^{\left|i-j\right|}\right)_{p\times p}=\begin{pmatrix}
 1 & r & \dots  & r^{p-1}\\
 r & 1 & \dots  & r^{p-2}\\
 \vdots  & \vdots  &  & \vdots \\
 r^{p-1} & r^{p-2} & \dots &1
\end{pmatrix},
    \end{equation*}
    where $r$ is fixed(i.e. $r$ do not change with $n$), and $0\le r <1$.\\

     $\mathbf{R}$ is the covariance matrix of the stationary $m$-dependent process, that is,\\
     \begin{equation*}
        \mathbf{R}=\left(r_{\left|i-j\right|}\right)_{p\times p}=\begin{pmatrix}
 1 & r_1 & \dots  & r_{p-1}\\
 r_1 & 1 & \dots  & r_{p-2}\\
 \vdots  & \vdots  &  & \vdots \\
 r_{p-1} & r_{p-2} & \dots &1
\end{pmatrix},
    \end{equation*}
where $\left \{ r_k \right \} ^\infty _{k=1}$ satisfies that there exist $m$ such that $r_k=0$ for all $k>m$.
\end{remark}

Let the distribution of $\xi$ be a Gumbel distribution with function
\begin{equation*}
    F_{\xi}\left ( x \right ) =\exp \left ( -\frac{1}{4\sqrt{2\pi } }e^{-x/2}  \right )  , x\in \mathbb{R}.
\end{equation*}
Define
\begin{equation*}
    \mu=2\sqrt{\log n}-\frac{\log \log n}{4\sqrt{\log n} } ,
\end{equation*}
and
\begin{equation*}
     a_p=8+\frac{16}{p} \sum_{k=1}^{p-1}r_{k}^2 \left ( p-k \right ).
\end{equation*}
First, we state he asymptotical distribution of the maximum interpoint Euclidean distance of high-dimensional normal distributions.
\begin{theorem}\label{th1}
    Under Condition (C1)-(C2) ,if $\bm{X}_i\overset{i.i.d.}{\sim} N\left(\bm \mu,\mathbf{R}\right)$, $i=1,\cdots,n$, then, as $n\rightarrow\infty$,
    \begin{equation*}
        4\sqrt{\log n}\left ( \frac{M_n^2-2p}{\sqrt{pa_p} } -\mu  \right )  \xrightarrow{d} \xi.
    \end{equation*}

\end{theorem}

In a more general case, we will consider the random vector $\bm{X}$ has dependent sub-gaussian components. Specifically, $\bm{X}_i=\mathbf{R}^{\frac{1}{2}}\bm{\epsilon}_i$, for any $i=1,2,\dots,n$, where $\bm{\epsilon}_i$'s are $p$-dimensional independently and identically distributed random vectors whose components are also independently and identically distributed with a sub-gaussian distribution. The definition of $\mathbf{R}$ is the same as that of formula (2).

Let $\mathbf{A}$ be a $n\times n$-dimension matrix whose element is write as $(\mathbf{A})_{ij}$ , $\left \| \mathbf{A} \right \|   :=\sum_{i=1}^{n} \sum_{j=1}^{n} \left | (\mathbf{A})_{ij} \right |$.
For two matrices $\mathbf{A}$ and $\mathbf{B}$ of the same dimension $m \times n$, the Hadamard product ${\displaystyle \mathbf{A}\odot \mathbf{B}}$  is a matrix of the same dimension as the operands, with elements given by
    \begin{equation*}
    \left ( \mathbf{A}\odot \mathbf{B} \right ) _{ij}=\left ( \mathbf{A}\right ) _{ij}\left ( \mathbf{B}\right ) _{ij}
    \end{equation*}
Furthermore, we write $\mathbf{A}^{\odot2}:=\mathbf{A}\odot \mathbf{A}$.
\begin{theorem}
  Assume $\bm{X}_i{=}\mathbf{R}^{\frac{1}{2}}\bm{\epsilon}_i$, for any $i=1,2,\dots,n$, where $\bm{X}_i$'s are $p$-dimensional i.i.d. random vectors. The components of $\bm{\epsilon}_i=(\epsilon_{i1},\cdots,\epsilon_{ip})$ are i.i.d. in a sub-gaussian distribution with $E(\epsilon_{ij})=0,var(\epsilon_{ij})=1,E(\epsilon_{ij}^4)=\kappa_4 < 5$.   Suppose Condition (C1)-(C2) hold.  Set $$b_p=p^{-1}\left ( 8\left \| \left ( \mathbf{T}^2 \right )^{\odot 2}  \right \| +2\left ( \kappa _4 -3\right ) \left \| \left ( \mathbf{T}^{\odot 2} \right )^2  \right \| \right ).$$ Then, we have, as $n\rightarrow\infty$,
    \begin{equation}\label{subgd}
        4\sqrt{\log n}\left ( \frac{M_n^2-2p}{\sqrt{pb_p}} -\mu  \right )  \xrightarrow{d} \xi ,
    \end{equation}
    where $\mathbf{T}:=\mathbf{R}^{\frac{1}{2}}$.
\end{theorem}
Note that if $\mathbf{R}=\mathbf{I}_p$, i.e. $r_k=0, k=1,\cdots,p-1$, we have $b_p=2(\kappa_4+1)$ and the limiting distribution (\ref{subgd}) is the same as (\ref{expd}).

We will list two applications of the result to high-dimensional data analysis.

1. {\it Test of covariance matrix.} Consider $\bm X_1,\cdots,\bm X_n\overset{i.i.d.}{\sim} N\left(\bm \mu,\mathbf{\Sigma}\right)$. We want to test the following hypothesis:
\begin{align}
H_0: \mathbf{\Sigma}=\mathbf{I}_p~v.s.~H_1:\mathbf{\Sigma}\not=\mathbf{I}_p.
\end{align}
We could adopt $M_n$ to test the above problem and reject the null hypothesis at the significant level $\alpha$ if
$$
M_n^2 \geq 2 p+\sqrt{p} \cdot \frac{q_\alpha+8 \log n-\log \log n}{\sqrt{2 \log n}}
$$
where $q_\alpha=-\log (32 \pi)-2 \log \log (1-\alpha)^{-1}$. In addition, if $\mathbf{\Sigma}=\mathbf{R}$ under the alternative hypothesis, the power function $\beta_{M_n}(r_n)$ is
\begin{align*}
\beta_{M_n}(r_n)=1-F_{\xi}\left(\sqrt{\frac{8}{a_p}}q_\alpha+\left(\sqrt{\frac{8}{a_p}}-1\right)(8\log n-\log\log n)\right).
\end{align*}
If $\frac{\log n}{p}\sum_{k=1}^{p-1}r_k^2(p-k) \to \infty$, we have $\beta_{M_n}(r_n)\to 1$. \citet{chen2010tests} proposed a test procedure based on a consistent estimator of ${\rm tr}[(\mathbf{\Sigma}-\mathbf{I}_p)^2]$ and show that it is consistent if $\frac{n}{p}{\rm tr}[(\mathbf{\Sigma}-\mathbf{I}_p)^2]\to \infty$. Note that if $\mathbf{\Sigma}=\mathbf{R}$,
${\rm tr}[(\mathbf{\Sigma}-\mathbf{I}_p)^2]=2\sum_{k=1}^{p-1}r_k^2(p-k)$. So the power of \citet{chen2010tests}'s test would converge to one if $\frac{n}{p}\sum_{k=1}^{p-1}r_k^2(p-k) \to \infty$, which indicates that \citet{chen2010tests}'s test is more efficient than $M_n$.

2. {\it outlier detection}. Suppose $\bm X_1,\cdots,\bm X_n\overset{i.i.d.}{\sim} N\left(\bm \mu,\mathbf{R}\right)$. We want to test the following hypothesis:
$$
\boldsymbol{H}_0 \text { : There is no typical outlier point among } \boldsymbol{X}_1, \boldsymbol{X}_2, \ldots, \boldsymbol{X}_n \text {. }
$$
According to Theorem \ref{th1}, we will reject the null hypothesis if
$$
M_n^2 \geq 2 p+\sqrt{\frac{pa_p}{8}} \cdot \frac{q_\alpha+8 \log n-\log \log n}{\sqrt{2 \log n}}.
$$
Note that $a_p=\frac{8}{p}{\rm tr(\mathbf{\Sigma}^2)}$ under the covariance assumption and $T_{2,n}$ is a consistent estimator of ${\rm tr(\mathbf{\Sigma}^2)}$ in \citet{chen2010tests}. So we can estimate $a_p$ by $\frac{8}{p}T_{2,n}$.

\section{Proof of main result}

Throughout the paper, $\bm{X}_1 , \dots , \bm{X}_n$ be identically distributed copies of $p$-dimensional random vector $\bm{X}$ with mean $\bm{\mu}$ and covariance matrix $\mathbf{R}$. Write $\bm{X}_i=\left(x_{i,1},x_{i,2},\dots,x_{i,p} \right)^T$ for each $i=1,2,\dots,n$, and $\bm{X}^{\left(j\right)}=\left(x_{1,j},x_{2,j},\dots,x_{n,j} \right)^T$ for each $j=1,2,\dots,p$. Write $y_{k,(i,j)}= x_{i,k}^2+x_{j,k}^2-2x_{i,k}x_{j,k} $ Thus,

\begin{equation*}
    M_n^2=\max_{1\le i<j\le n} \sum_{k=1}^{p}\left ( x_{i,k}^2+x_{j,k}^2-2x_{i,k}x_{j,k} \right )
    =\max_{1\le i<j\le n} \sum_{k=1}^{p} y_{k,(i,j)}.
\end{equation*}

We state the following three important lemmas here.

\begin{lemma} [\citet{ref9}, Theorem 1] Let $\bm{Y}_1 ,\dots, \bm{Y}_n$ be a sequence of p-dimensional dependent random vectors with mean zero, i.e., $\mathbb{E}\left(\bm{Y}_t\right)=0$ Write $\bm{S}_n=n^{-1/2} {\textstyle \sum_{t=1}^{n}}\bm{Y}_t $ and $\mathbf{\Xi} =\mathrm{Cov} \left ( S_n \right ) $. \\

Condition a: There exist a sequence of constants $B_n\ge 1$ and a universal constant $\gamma_1 \ge 1$ such that $\left \| Y_{t,j} \right \|  _{\psi _{\gamma _1}}\le B_n $ for all $t \in \left \{ 1,2,\dots , n \right \} $ and $j \in \left \{ 1,2,\dots , p \right \} $. \\

Condition b: There exist some universal constants $K_1>1$ , $K_2>0$ and $\gamma_2>0$ such that $\alpha _n\left ( k \right ) \le K_1\exp \left ( -K_2k^{\gamma _2} \right ) $ for any $k\ge 1$.\\

Condition c: There exist a universal constant $K_3$ such that $\min \text{diag}\left(\mathbf{\Xi}\right)\ge K_3$.\\

Assume $\{\bm{Y}_t\}$ is an $\alpha$-mixing sequence with $p\ge n^\kappa$ for some universal constant $\kappa >0$. Under Condition a,b,c, it holds that
    \begin{equation*}
    \sup_{A\in \mathcal{A}}\left | \mathbb{P}\left ( \bm{S}_n\in A \right ) - \mathbb{P}\left ( \bm{G}\in A \right )  \right |\lesssim \frac{B_n^{2/3}\left ( \log p \right )^{\left ( 1+2\gamma _2 \right ) /\left ( 3\gamma _2 \right ) } }{n^{1/9}} +\frac{B_n\left ( \log p \right )^{7/6 } }{n^{1/9}},
    \end{equation*}
where $G\sim N\left(0,\Xi \right)$ and $\mathcal{A}$ is a class of Borel subsets in $\mathbb{R}^p$.\\

\end{lemma}

\begin{lemma}[\citet{ref10}, Theorem 1]
 Let  $\left\{\eta_\alpha,\alpha\in I\right\}$ be random variables on an index set $I$ and $\left\{B_\alpha,\alpha\in I\right\}$ be a set of subsets of $I$ , that is, for each $\alpha\in I$ , $B_\alpha\subset I$. for any $t\in \mathbb{R}$ , set $\lambda_p=\sum_{\alpha\in I}\mathbb{P}\left(\eta_\alpha>t\right)$, Then we have
 \begin{equation*}
     \left | \mathbb{P}\left (\max_{\alpha \in I}\eta _\alpha \le t \right )-e^{-\lambda_p}  \right | \le\left ( 1\wedge \lambda _p^{-1} \right ) \left ( b_1+b_2+b_3 \right ) ,
 \end{equation*}
 where
 \begin{align*}
     b_1&=\sum_{\alpha \in I}\sum_{\beta\in B_\alpha } \mathbb{P}\left (  \eta_\alpha>s_n \right )\mathbb{P}\left (  \eta_\beta>s_n \right ) ,\\
     b_2&=\sum_{\alpha \in I}\sum_{\alpha\ne\beta\in B_\alpha } \mathbb{P}\left (  \eta_\alpha>s_n ,  \eta_\beta>s_n \right ),\\
     b_3&=\sum_{\alpha \in I}\left | \mathbb{P}\left \{ \eta_\alpha >t|\sigma \left ( \eta_\beta ,\beta \notin B_\alpha  \right )  \right \}-\mathbb{P}\left ( \eta_\alpha >t \right )   \right | ,
 \end{align*}
 and $\sigma \left ( \eta_\beta ,\beta \notin B_\alpha  \right ) $ is the $\sigma$-algebra generated by $ \left ( \eta_\beta ,\beta \notin B_\alpha  \right )$. In particular, if $\eta_\alpha$ is independent of $ \left ( \eta_\beta ,\beta \notin B_\alpha  \right )$ for each $\alpha$, then $b_3$ vanishes.
 \end{lemma}

\begin{lemma}
Let $A_n$ be a random variable for each $n\ge1$ satisfying
\begin{equation*}
    \lim_{n \to \infty} \mathbb{P}\left ( A_n\le \sqrt{4\log n-\log \log n+x} \right ) = F\left(x\right)
\end{equation*}
for any $x\in \mathbb{R}$, where $F\left(x\right)$ is a continuous distribution function on $\mathbb{R}$. Then\\
\begin{equation*}
    A_n=2\sqrt{\log n}-\frac{\log \log n}{4\sqrt{\log n} }+  \frac{U_n}{4\sqrt{\log n} },
\end{equation*}
where $U_n$ converges weakly to a probability measure with distribution function $F\left(x\right)$.
\end{lemma}

\subsection{Proof of Theorem 2.1}

Step 1: find the covariance structure of $y_{k,\left ( i,j \right ) }$.\\Obviously $\mathbb{E}\left ( y_{k,\left ( i,j \right ) } \right ) =2$ , $\mathrm{Var}\left ( y_{k,\left ( i,j \right ) } \right )=8 $ and
\begin{align*}
  \mathrm{Cov}\left (  y_{k,\left ( i,j \right )}, y_{m,\left ( i,j \right )} \right )
  =&  \mathrm{Cov}\left ( x_{i,k}^2+x_{j,k}^2-2x_{i,k}x_{j,k},x_{i,m}^2+x_{j,m}^2-2x_{i,m}x_{j,m} \right ) \\
  =&\mathrm{Cov}\left ( x_{i,k}^2,x_{i,m}^2 \right )+\mathrm{Cov}\left ( x_{j,k}^2,x_{j,m}^2 \right )+4\mathrm{Cov}\left ( x_{i,k}x_{j,k},x_{i,m}x_{j,m} \right )\\
  =&2\mathbb{E}\left ( x_{i,k}^2x_{i,m}^2 \right )-2+4\mathbb{E}\left ( x_{i,k}x_{j,k}x_{i,m}x_{j,m} \right )\\
  =&2\mathbb{E}\left ( x_{i,k}^2x_{i,m}^2 \right )-2+4r^2_{\left | k-m \right | }.
\end{align*}
Let $\eta_1, \eta_2\sim N\left(0,1\right)$. Set $a_1,a_2$ satisfy $a_1^2+a_2^2=1$ and $2a_1a_2=r_{\left | k-m \right | }$ we have $\left (   x_{i,m}, x_{i,k}\right ) \overset{d}{=}\left ( a_1\eta_1+a_2\eta_2, a_2\eta_1+a_1\eta_2\right ) $
Thus,
\begin{align*}
    \mathbb{E}\left ( x_{i,k}^2x_{i,m}^2 \right )&=\mathbb{E}\left ( a_1^2a_2^2\left ( \eta _1^4+\eta_2^4 \right )+\left ( a_1^4+4a_1^2a_2^2+a_2^4 \right )\eta^2_1\eta_2^2   \right )\\
    &=2r^2_{\left | k-m \right | }+1.
\end{align*}
So $\mathrm{Cov}\left (  y_{k,\left ( i,j \right )}, y_{m,\left ( i,j \right )} \right )=8r^2_{\left | k-m \right | }$,
\begin{align*}
   \mathrm{Cov}\left (  y_{k,\left ( i,j \right )}, y_{k,\left ( i,l \right )} \right )
  =&  \mathrm{Cov}\left ( x_{i,k}^2+x_{j,k}^2-2x_{i,k}x_{j,k},x_{i,k}^2+x_{l,k}^2-2x_{i,k}x_{l,k} \right ) \\
  =&\mathrm{Cov}\left ( x_{i,k}^2,x_{i,k}^2 \right ) =2
\end{align*}
and
\begin{align*}
   \mathrm{Cov}\left (  y_{k,\left ( i,j \right )}, y_{m,\left ( i,l \right )} \right )
  =&  \mathrm{Cov}\left ( x_{i,k}^2+x_{j,k}^2-2x_{i,k}x_{j,k},x_{i,m}^2+x_{l,m}^2-2x_{i,m}x_{l,m} \right ) \\
  =&\mathrm{Cov}\left ( x_{i,k}^2,x_{i,m}^2 \right ) =2r^2_{\left | k-m \right | }.
\end{align*}

Let $\bm{Y}_t=\left(Y_{t,1},Y_{t,2},\dots ,Y_{t,\frac{n\left ( n-1 \right )}{2}}\right)^T$ be the $\frac{n\left ( n-1 \right )}{2} $-dimensional random straightened vector constructed from $y_{t,\left ( i,j \right )}$ ,$1\le i<j\le n$. The subscript $\left(i,j\right)$ is converted into subscript $l$ by straightening. Denote this subscript mapping as $h:\left(i,j\right)\mapsto l$ . Obviously it is a bijective mapping. Write $\left ( \sigma _{l,k} \right ) _{\frac{n\left ( n-1 \right )}{2}\times \frac{n\left ( n-1 \right )}{2}}=\mathbf{\Sigma}:=\mathrm{Cov} \left ( \bm{Y}_1 \right ) /8$,  Through the previous discussion, we can calculate $\sigma _{l,k}$'s satisfy the following formula.
\begin{equation*}
    \sigma _{l,k}=\begin{cases}
 1 & \text{ if } l=k \\
 1/4 & \text{ if } \exists i\text{ appear in both } h^{-1}\left ( l \right ) \text{ and }  h^{-1}\left ( k \right ) \\
 0 & \text{ others }
\end{cases}
\end{equation*}
With the function $h$ held constant, $\mathbf{\Sigma}$ is exclusively associated with $n$. Then, write $\bm{S}_p=p^{-1/2} {\textstyle \sum_{t=1}^{p}}\bm{Y}_t $ ,Thus
\begin{align*}
    \mathbf{\Xi}:=&\mathrm{Cov}\left(\bm{S}_p\right)\\
    =&p^{-1}\left ( \sum_{t=1}^{p} \mathrm{Cov}\left ( \bm{Y}_t \right ) +2\sum_{1\le t<k\le p} \mathrm{Cov}\left ( \bm{Y}_t,\bm{Y}_k \right )    \right )\\
    =&p^{-1}\left ( p\mathrm{Cov}\left ( \bm{Y}_1 \right ) +2\sum_{1\le t<k\le p} r^2_{\left|t-k\right|}\mathrm{Cov}\left ( \bm{Y}_1 \right )    \right )\\
    =&p^{-1}\left ( p+2\sum_{i=1}^{p-1}r^2_{i} \left ( p-i \right )  \right ) \mathrm{Cov}\left ( \bm{Y}_1 \right )\\
    =&a_p\mathbf{\Sigma},
\end{align*}
where $a_p$ defined in (5).

Step 2:
Check Condition a,b,c in Lemma 1.

Notice that $Y_{t,j}/2\sim \chi^2\left(1\right)$ for all $t \in \left \{ 1,2,\dots , p\right \} $ and $j \in \left \{ 1,2,\dots , \frac{n\left(n+1\right)}{2} \right \} $. Therefore $Y_{t,j}$'s are non-negative sub-exponential distribution random variables. So there exist universal constant $B$ such that $\left \| Y_{t,j} \right \|  _{\psi _{\gamma _1}}\le B$ Thus condition a in Lemma 1 holds.

By Condition (C1), we have
\begin{equation*}
    \sup_{t}\sup_{f\in L^2\left(\mathcal{F}_{-\infty}^t\right), g\in L^2\left(\mathcal{F}^{\infty}_{t+k}\right) }\left|\mathrm{corr}\left(f,g\right)\right| \le K_1\exp \left ( -K_2k^{\gamma } \right ) ,
\end{equation*}
where $\mathcal{F}_{-\infty}^u $ and $\mathcal{F}^{\infty}_u$ the $\sigma$-fields generated respectively by $\left\{\bm{X}^{\left(t\right)}\right\}_{t\le u}$ and $\left\{\bm{X}^{\left(t\right)}\right\}_{t\ge u}$. By the definition of $\bm{Y}_t$ we know $\bm{Y}_t\in \sigma\left(\bm{X}^{\left(t\right)}\right)$ . Write the $\alpha$-mixing coefficient and $\rho$-mixing coefficient of $\bm{Y}_t$ as $\alpha_n'$ and $\rho_n'$. Write $\mathcal{G}_{-\infty}^u $ and $\mathcal{G}^{\infty}_u$ the $\sigma$-fields generated respectively by $\left\{\bm{Y}^{\left(t\right)}\right\}_{t\le u}$ and $\left\{\bm{Y}^{\left(t\right)}\right\}_{t\ge u}$. So we have $L^2\left(\mathcal{G}_{-\infty}^u \right)\subset L^2\left(\mathcal{F}_{-\infty}^u \right)$ and $L^2\left(\mathcal{G}^{\infty}_u \right)\subset L^2\left(\mathcal{F}^{\infty}_u \right)$ . Thus,
\begin{equation*}
    \sup_{t}\sup_{f\in L^2\left(\mathcal{G}_{-\infty}^t\right), g\in L^2\left(\mathcal{G}^{\infty}_{t+k}\right) }\left|\mathrm{corr}\left(f,g\right)\right| \le \sup_{t}\sup_{f\in L^2\left(\mathcal{F}_{-\infty}^t\right), g\in L^2\left(\mathcal{F}^{\infty}_{t+k}\right) }\left|\mathrm{corr}\left(f,g\right)\right|.
\end{equation*}
Therefore, $\alpha'_n\left(k\right)\le\rho'_n\left(k\right)\le\rho_n\left(k\right)\le K_1\exp \left ( -K_2k^{\gamma } \right )$. So the condition b in Lemma 1 holds.

$\min \text{diag}\left(\mathbf{\Xi}\right)=a_p=8+\frac{16}{p} \sum_{k=1}^{p-1}r_{k}^2 \left ( p-k \right )>8$ So the condition c in Lemma 1 holds. By Lemma 1
\begin{equation*}
    \sup_{A\in \mathcal{A}}\left | \mathbb{P}\left ( \bm{S}_p-2\bm{1}_{\frac{n\left(n-1\right)}{2}} \in A \right ) - \mathbb{P}\left ( \bm{G}\in A \right )  \right |\lesssim \frac{B^{2/3}\log\left ( \frac{n\left ( n-1 \right ) }{2}  \right )^{\frac{1+2\gamma }{3\gamma } } }{p^{\frac{1}{9} }}+\frac{B\log\left ( \frac{n\left ( n-1 \right ) }{2}  \right )^{\frac{7}{6} } }{p^{\frac{1}{9} }},
\end{equation*}
where $\bm{G}\sim N\left(0,\mathbf{\Xi} \right)$ and $\mathcal{A}$ is a class of Borel subsets in $\mathbb{R}^{\frac{n\left(n-1\right)}{2}}$. By calculating the right formula, we can find $\mathrm{RHS}\asymp \frac{\left ( \log n \right )^{\nu} }{p^{\frac{1}{9} }} $ where $\nu=\frac{2+4\gamma }{3\gamma }\vee \frac{7}{3}$. By Condition (C2), we have $\mathrm{RHS}=o\left(1\right)$.\\
Notice that $M^2_n=\sqrt{p}\max_{1\le i<j\le n}\bm{S}_p$, Thus\\
\begin{equation*}
    \sup_{t}\left | \mathbb{P}\left ( M^2_n-2p \le t \right ) - \mathbb{P}\left ( \max_{1\le i<j\le n}\bm{G} \le t \right )  \right |=o\left(1\right).
\end{equation*}
Let $t=\sqrt{a_p}x$ , $A_n=\max_{1\le i<j\le n}{a_p}^{-1/2}\bm{G}$ , Thus
\begin{equation*}
    \sup_{x}\left | \mathbb{P}\left ( \frac{M^2_n-2p}{\sqrt{pa_p}} \le x \right ) - \mathbb{P}\left ( A_n \le x \right )  \right |=o\left(1\right).
\end{equation*}
Step 3: \\
Set $I=\left \{ \left ( i,j \right ) ;1\le i <j\le n \right \} $ for $\alpha= \left ( i,j \right ) \in I$
define $$B_\alpha=\left \{ \left ( k,l \right )\in I ;\left \{ k,l \right \}\cap \left \{ i,j \right \}\ne\emptyset ,\left ( k,l \right )\ne\alpha  \right \}$$ Let $\eta_{\left ( i,j \right )}$ be the term corresponding to $\left ( i,j \right )$ in ${a_p}^{-1/2}G$. We can find that $\eta_{\left ( i,j \right )}\sim N\left(0,1\right)$ and the covariance structure of it is $\mathbf{\Sigma}$.\\
Set $s_n=\sqrt{4\log n-\log \log n+x}$ By Lemma 2, we have\\
\begin{equation*}
\left | \mathbb{P}\left ( \max_{\alpha \in I}\eta_\alpha \le s_n \right ) -e^{-\lambda _p} \right | \le u_1+u_2,
\end{equation*}
where
\begin{align*}
    \lambda_p&=\sum_{\alpha \in I}\mathbb{P}\left (  \eta_\alpha>s_n \right ) =\frac{n\left ( n-1 \right ) }{2} \left ( 1-\Phi \left ( s_n \right )  \right ) , \\
    u_1&=\sum_{\alpha \in I}\sum_{\beta\in B_\alpha } \mathbb{P}\left (  \eta_\alpha>s_n \right )\mathbb{P}\left (  \eta_\beta>s_n \right ) =\frac{n\left ( n-1 \right ) }{2}\cdot \left ( 2\left ( n-2 \right )  \right )  \left ( 1-\Phi \left ( s_n \right )  \right ) ^2 ,\\
    u_2&=\sum_{\alpha \in I}\sum_{\beta\in B_\alpha } \mathbb{P}\left (  \eta_\alpha>s_n ,  \eta_\beta>s_n \right ) =\frac{n\left ( n-1 \right ) }{2}\cdot \left ( 2\left ( n-2 \right )  \right ) \cdot  \mathbb{P}\left (  \eta_\alpha>s_n ,  \eta_\beta>s_n \right ).
\end{align*}
By the formula $1-\Phi \left ( x \right )=\frac{1}{\sqrt{2\pi }x } e^{-x^2/2}\left ( 1+o\left ( 1 \right )  \right )$ as $x\rightarrow \infty$ , we conclude as $n\rightarrow\infty$ \\
\begin{align*}
    \lambda_p&=\frac{n^2}{2}\cdot  \left ( 1-\Phi \left ( s_n \right )  \right ) \left ( 1+o\left ( 1 \right )  \right )
    =\frac{n^2}{2\sqrt{2\pi }s_n } e^{-s_n^2/2}\left ( 1+o\left ( 1 \right )  \right )
    \rightarrow \frac{1}{4\sqrt{2\pi } } e^{x/2}.
\end{align*}
We can see
\begin{equation}
    \left ( 1-\Phi \left ( s_n \right )  \right )=O_p\left(n^{-2}\right).
\end{equation}
Thus $u_1=O\left(n^{-1}\right)\rightarrow0$ as $n\rightarrow\infty$ .\\
Then we estimate $u_2$ . Notice that $\left ( \eta_\alpha,\eta_\beta  \right )^T \sim N\left ( 0,\begin{pmatrix}
  1& \frac{1}{4} \\
 \frac{1}{4}  &1
\end{pmatrix} \right ) $ we can calculate\\
\begin{align*}
    u_2&\le\frac{n\left ( n-1 \right ) }{2}\cdot \left ( 2\left ( n-2 \right )  \right ) \cdot  \mathbb{P}\left (  \eta_\alpha+  \eta_\beta>2s_n \right )\\
    &\le\frac{n\left ( n-1 \right ) }{2}\cdot \left ( 2\left ( n-2 \right )  \right ) \cdot  \mathbb{P}\left (  N\left ( 0,\frac{5}{2}  \right ) >2s_n \right )\\
    &=n^3\cdot\left ( 1-\Phi \left ( \frac{2\sqrt{10}s_n }{5}  \right )  \right ) \cdot\left ( 1+o\left ( 1 \right )  \right )\\
    & \asymp \frac{n^3}{s_n} e^{-\frac{4s_n^2}{5} }=o\left ( n^{-1/5} \right ) \rightarrow 0.
\end{align*}
Thus, for any $x\in \mathbb{R}$
\begin{equation*}
    \lim_{n \to \infty} \mathbb{P}\left ( A_n\le s_n \right ) =\lim_{n \to \infty} \mathbb{P}\left ( \max_{\alpha \in I}\eta_\alpha \le s_n \right ) =\exp \left ( -\frac{1}{4\sqrt{2\pi } }e^{-x/2}  \right ).
\end{equation*}
Then, by Lemma 3,
\begin{equation*}
    A_n=2\sqrt{\log n}-\frac{\log \log n}{4\sqrt{\log n} }+  \frac{U_n}{4\sqrt{\log n} } ,
\end{equation*}
where $U_n\overset{d}{\rightarrow}\xi$ with distribution function $F_{\xi}\left ( x \right )$ in (3.1). So we have\\
\begin{align*}
    M_n^2=2p+\sqrt{pa_p}A_n=2p+\sqrt{pa_p}\left(2\sqrt{\log n}-\frac{\log \log n}{4\sqrt{\log n} }+  \frac{U_n}{4\sqrt{\log n} }\right).
\end{align*}
It follows that
\begin{align*}
    4\sqrt{\log n}\left ( \frac{M_n^2-2p}{\sqrt{pa_p} } -2\sqrt{\log n}-\frac{\log \log n}{4\sqrt{\log n} }  \right )  \xrightarrow{d} \xi.
\end{align*}

\subsection{Proof of Theorem 2.2}

Step 1: find the covariance structure of $y_{k,\left(i,j\right)}$\\
Similarly in proof of theorem 2.1, we can get
\begin{align*}
   \mathbb{E}\left(y_{k,\left(i,j\right)}\right)&=2\\
   \mathrm{Cov}\left(y_{k,\left(i,j\right)},y_{k,\left(i,j\right)}\right)&=2\mathbb{E}\left(x_{i,k}^4\right)+2\\
   \mathrm{Cov}\left(y_{k,\left(i,j\right)},y_{k,\left(i,l\right)}\right)&=\mathbb{E}\left(x_{i,k}^4\right)-1\\
   \mathrm{Cov}\left(y_{k,\left(i,j\right)},y_{m,\left(i,j\right)}\right)&=2\mathbb{E}\left(x_{i,k}^2x_{i,m}^2\right)+4r^2_{\left|k-m\right|}-2\\
   \mathrm{Cov}\left(y_{k,\left(i,j\right)},y_{m,\left(i,l\right)}\right)&=\mathbb{E}\left(x_{i,k}^2x_{i,m}^2\right)-1,
\end{align*}
where $k\ne m$ and $i,j,l$ are not equal to each other. Notice that $\bm{X}_i\overset{d}{=}\mathbf{T}\bm{\epsilon}$. Write the $k^{th}$ row of $\mathbf{T}$ as $\bm{T}_k=\left(t_{k1},t_{k2},\dots,t_{kp}\right)$ , thus
\begin{equation}
x_{i,k}\overset{d}{=}\bm{T}_k\bm{\epsilon}=\sum_{j=1}^{p} t_{kj}\epsilon_j .
\end{equation}
Similarly we can define $x_{i,m}$. Then
\begin{align*}
    \mathbb{E}\left(x^4_{i,k}\right)&=\mathbb{E}\left(\left(\sum_{j=1}^{p} t_{kj}\epsilon_j\right)^4\right)\\
    &=\mathbb{E}\left(\left(\sum_{j=1}^{p} t^2_{kj}\epsilon^2_j+2\sum_{1\le j < l \le p}t_{kj}t_{kl}\epsilon_j\epsilon_l\right)^2\right)\\
    &=\mathbb{E}\left(\left(\sum_{j=1}^{p} t^2_{kj}\epsilon^2_j\right)^2+\left(2\sum_{1\le j < l \le p}t_{kj}t_{kl}\epsilon_j\epsilon_l\right)^2\right)\\
    &=\mathbb{E}\left(\left(\sum_{j=1}^{p}t_{kj}^4\epsilon_j^4\right)+\left(2\sum_{1\le j<l\le p}t^2_{kj}t^2_{kl}\epsilon_j^2\epsilon_l^2\right)+\left(4\sum_{1\le j<l\le p}t^2_{kj}t^2_{kl}\epsilon_j^2\epsilon_l^2\right)\right)\\
    &=\left(\sum_{j=1}^{p}t_{kj}^4\right)\kappa_4+6\sum_{1\le j<l\le p}t^2_{kj}t^2_{kl}
\end{align*}
and
\begin{align*}
    \mathbb{E}\left(x^2_{i,k}x^2_{i,m}\right)&=\mathbb{E}\left(\left(\sum_{j=1}^{p} t^2_{kj}\epsilon^2_j+2\sum_{1\le j < l \le p}t_{kj}t_{kl}\epsilon_j\epsilon_l\right)\left(\sum_{j=1}^{p} t^2_{mj}\epsilon^2_j+2\sum_{1\le j < l \le p}t_{mj}t_{ml}\epsilon_j\epsilon_l\right)\right)\\
    &=\mathbb{E}\left(\sum_{j=1}^{p}t^2_{kj}t^2_{mj}\epsilon^4+\sum_{1\le j,l\le p,j\ne l} t^2_{kj}t^2_{ml}\epsilon_j^2\epsilon_l^2+ 4\sum_{1\le j < l \le p}t_{kj}t_{kl}t_{mj}t_{ml}\epsilon_j^2\epsilon^2_l\right)\\
    &=\left(\sum_{j=1}^{p}t^2_{kj}t^2_{mj}\right)\kappa_4+\sum_{1\le j,l\le p,j\ne l} t^2_{kj}t^2_{ml}+ 4\sum_{1\le j < l \le p}t_{kj}t_{kl}t_{mj}t_{ml}  .
\end{align*}
Notice $\mathbf{T}^2=\mathbf{R}$, we have $\sum_{j=1}^{p} t^2_{kj}=\sum_{j=1}^{p} t^2_{mj}=1$, $\sum_{j=1}^{p} t_{kj}t_{mj}=r_{\left | k-m \right | }$. Thus
\begin{align*}
    2\sum_{1\le j<l\le p}t^2_{kj}t^2_{kl}&=\left(\sum_{j=1}^{p} t^2_{kj}\right)^2-\sum_{j=1}^{p}t_{kj}^4=1-\sum_{j=1}^{p}t_{kj}^4\\
    \sum_{1\le j,l\le p,j\ne l} t^2_{kj}t^2_{ml}&=\left(\sum_{j=1}^{p} t^2_{kj}\right)\left(\sum_{j=1}^{p} t^2_{mj}\right)-\sum_{j=1}^{p}t^2_{kj}t^2_{mj}=1-\sum_{j=1}^{p}t^2_{kj}t^2_{mj}\\
    2\sum_{1\le j < l \le p}t_{kj}t_{kl}t_{mj}t_{ml}&=\left(\sum_{j=1}^{p} t_{kj}t_{mj}\right)^2-\sum_{j=1}^{p}t^2_{kj}t^2_{mj}=r^2_{\left | k-m \right | }-\sum_{j=1}^{p}t^2_{kj}t^2_{mj}.
\end{align*}
So $\mathbb{E}\left(x^4_{i,k}\right)=\left(\sum_{j=1}^{p}t_{kj}^4\right)\left(\kappa_4-3\right)+3$ , $\mathbb{E}\left(x^2_{i,k}x^2_{i,m}\right)=\left(\sum_{j=1}^{p}t^2_{kj}t^2_{mj}\right)\left(\kappa_4-3\right)+2r^2_{\left | k-m \right | }+1$.
Write $(c_{ij})_{p\times p}=\mathbf{C}=\left ( \mathbf{T}^{\odot 2} \right )^2 $. We have $c_{kk}=\sum_{j=1}^{p}t_{kj}^4$ and  $c_{km}=\sum_{j=1}^{p}t^2_{kj}t^2_{mj}$
Thus,
\begin{align*}
   \mathrm{Cov}\left(y_{k,\left(i,j\right)},y_{k,\left(i,j\right)}\right)&=8+2c_{kk}\left(\kappa_4-3\right)\\
   \mathrm{Cov}\left(y_{k,\left(i,j\right)},y_{k,\left(i,l\right)}\right)&=2+c_{kk}\left(\kappa_4-3\right)\\
   \mathrm{Cov}\left(y_{k,\left(i,j\right)},y_{m,\left(i,j\right)}\right)&=8r^2_{\left | k-m \right | }+2c_{km}\left(\kappa_4-3\right)\\
   \mathrm{Cov}\left(y_{k,\left(i,j\right)},y_{m,\left(i,l\right)}\right)&=2r^2_{\left | k-m \right | }+c_{km}\left(\kappa_4-3\right).
\end{align*}
Let $\bm{Y}_t=\left(Y_{t,1},Y_{t,2},\dots ,Y_{t,\frac{n\left ( n-1 \right )}{2}}\right)^T$ be the $\frac{n\left ( n-1 \right )}{2} $-dimensional random straightened vector constructed from $y_{t,\left ( i,j \right )}$ ,$1\le i<j\le n$. The subscript $\left(i,j\right)$ is converted into subscript $l$ by straightening. Denote this subscript mapping as $h:\left(i,j\right)\mapsto l$ . Obviously it is a bijective mapping. Write $\mathbf{\Sigma}_1:=\left ( \sigma _{1,lk} \right ) _{\frac{n\left ( n-1 \right )}{2}\times \frac{n\left ( n-1 \right )}{2}}$
\begin{equation*}
    \sigma _{1,lk}=\begin{cases}
 1 & \text{ if } l=k \\
 1/4 & \text{ if } \exists i\text{ appear in both } h^{-1}\left ( l \right ) \text{ and }  h^{-1}\left ( k \right ) \\
 0 & \text{ others }
\end{cases}
\end{equation*}
and $\mathbf{\Sigma}_2:=\left ( \sigma _{2,lk} \right ) _{\frac{n\left ( n-1 \right )}{2}\times \frac{n\left ( n-1 \right )}{2}}$
\begin{equation*}
    \sigma _{2,lk}=\begin{cases}
 1 & \text{ if } l=k \\
 1/2 & \text{ if } \exists i\text{ appear in both } h^{-1}\left ( l \right ) \text{ and }  h^{-1}\left ( k \right ) \\
 0 & \text{ others }
\end{cases}.
\end{equation*}

With the function $h$ held constant, $\mathbf{\Sigma}_1$ and $\mathbf{\Sigma}_2$ are exclusively associated with $n$. Through the previous discussion, we can calculate $\mathrm{Cov}\left(\bm{Y}_t\right)=8\mathbf{\Sigma}_1+2c_{tt}\left(\kappa_4-3\right)\mathbf{\Sigma}_2$ and $\mathrm{Cov}\left(\bm{Y}_t,\bm{Y}_k\right)=8r^2_{\left|t-k\right|}\mathbf{\Sigma}_1+2c_{tk}\left(\kappa_4-3\right)\mathbf{\Sigma}_2$ Then, write $\bm{S}_p=p^{-1/2} {\textstyle \sum_{t=1}^{p}}\bm{Y}_t $ ,Thus
\begin{align*}
    \mathbf{\Xi}:=&\mathrm{Cov}\left(\bm{S}_p\right)\\
    =&p^{-1}\left ( \sum_{t=1}^{p} \mathrm{Cov}\left ( \bm{Y}_t \right ) +2\sum_{1\le t<k\le p} \mathrm{Cov}\left ( \bm{Y}_t,\bm{Y}_k \right )    \right )\\
    =&p^{-1}\left ( \sum_{t=1}^{p} 8\mathbf{\Sigma}_1+2c_{tt}\left(\kappa_4-3\right)\mathbf{\Sigma}_2 +2\sum_{1\le t<k\le p} 8r^2_{\left|t-k\right|}\mathbf{\Sigma}_1+2c_{tk}\left(\kappa_4-3\right)\mathbf{\Sigma}_2\right )\\
    =&p^{-1}\left ( 8\left \| \left ( \mathbf{T}^2 \right )^{\odot 2}  \right \|\mathbf{\Sigma}_1 +2\left ( \kappa _4 -3\right ) \left \| \left ( \mathbf{T}^{\odot 2} \right )^2  \right \| \mathbf{\Sigma}_2\right ) \\
    =&b_p\mathbf{\Sigma}_p,
\end{align*}
where $b_p$ defined as $p^{-1}\left ( 8\left \| \left ( \mathbf{T}^2 \right )^{\odot 2}  \right \| +2\left ( \kappa _4 -3\right ) \left \| \left ( \mathbf{T}^{\odot 2} \right )^2  \right \| \right )$ , $\mathbf{\Sigma}_p:=\left ( \sigma _{p,lk} \right ) _{\frac{n\left ( n-1 \right )}{2}\times \frac{n\left ( n-1 \right )}{2}}$
\begin{equation*}
    \sigma _{p,lk}=\begin{cases}
 1 & \text{ if } l=k \\
 \rho_p & \text{ if } \exists i\text{ appear in both } h^{-1}\left ( l \right ) \text{ and }  h^{-1}\left ( k \right ) \\
 0 & \text{ others }
\end{cases},
\end{equation*}
where
\begin{equation}
    \rho_p=\frac{2\left \| \left ( \mathbf{T}^2 \right )^{\odot 2}  \right \| +1\left ( \kappa _4 -3\right ) \left \| \left ( \mathbf{T}^{\odot 2} \right )^2  \right \|}{8\left \| \left ( \mathbf{T}^2 \right )^{\odot 2}  \right \| +2\left ( \kappa _4 -3\right ) \left \| \left ( \mathbf{T}^{\odot 2} \right )^2  \right \|}.
\end{equation}
Notice that $\left \| \left ( \mathbf{T}^2 \right )^{\odot 2}  \right \|=\sum_{1\le i,j\le p}\left(\sum_{k=1}^{p}t_{ik}t_{kj}\right)^2$ and $\left \| \left ( \mathbf{T}^{\odot 2} \right )^2  \right \|=\sum_{1\le i,j\le p}\left(\sum_{k=1}^{p}t^2_{ik}t^2_{kj}\right)$. Obviously $\left \| \left ( \mathbf{T}^2 \right )^{\odot 2}  \right \|\ge \left \| \left ( \mathbf{T}^{\odot 2} \right )^2  \right \|$. By assuming $\kappa\le 5$, we can get $\rho_p\le \frac{1}{3}$.\\
Step 2:
Check Condition a,b,c in Lemma 1

Notice that $Y_{t,j}=\left(x_{k,t}-x_{m,t}\right)^2$ for some $m$ and $k$. Recall formula (8), we have
\begin{equation*}
    x_{k,t}-x_{m,t}\overset{d}{=}\sum_{j=1}^{p} t_{tj}\epsilon^k_j-\sum_{j=1}^{p} t_{tj}\epsilon^m_j,
\end{equation*}
where $\epsilon^k_j$ and $\epsilon^m_j$ are i.i.d. copy of $\epsilon_j$. Since $\epsilon_j$'s are sub-gaussian with parameter $\sigma$, $\epsilon^k_j-\epsilon^m_j$'s are also sub-gaussian with same parameter $\sqrt{2}\sigma$. and $x_{k,t}-x_{m,t}$ is a linear combination of $\epsilon^k_j-\epsilon^m_j$,  $x_{k,t}-x_{m,t}$ is sub-gaussian with parameter $\sqrt{2\sum_{j=1}^p t^2_{tj}}\sigma=\sqrt{2}\sigma$
Therefore $Y_{t,j}$'s are non-negative sub-exponential distribution random variables. So there exist universal constant $B$ such that $\left \| Y_{t,j} \right \|  _{\psi _{\gamma _1}}\le B$ Thus condition a in Lemma 1 holds.

By Condition (C1), we have
\begin{equation*}
    \sup_{t}\sup_{f\in L^2\left(\mathcal{F}_{-\infty}^t\right), g\in L^2\left(\mathcal{F}^{\infty}_{t+k}\right) }\left|\mathrm{corr}\left(f,g\right)\right| \le K_1\exp \left ( -K_2k^{\gamma } \right ) ,
\end{equation*}
where $\mathcal{F}_{-\infty}^u $ and $\mathcal{F}^{\infty}_u$ the $\sigma$-fields generated respectively by $\left\{\bm{X}^{\left(t\right)}\right\}_{t\le u}$ and $\left\{\bm{X}^{\left(t\right)}\right\}_{t\ge u}$. By the definition of $Y_t$ we know $\bm{Y}_t\in \sigma\left(\bm{X}^{\left(t\right)}\right)$ . Write the $\alpha$-mixing coefficient and $\rho$-mixing coefficient of $\bm{Y}_t$ as $\alpha_n'$ and $\rho_n'$. Write $\mathcal{G}_{-\infty}^u $ and $\mathcal{G}^{\infty}_u$ the $\sigma$-fields generated respectively by $\left\{\bm{Y}^{\left(t\right)}\right\}_{t\le u}$ and $\left\{\bm{Y}^{\left(t\right)}\right\}_{t\ge u}$. So we have $L^2\left(\mathcal{G}_{-\infty}^u \right)\subset L^2\left(\mathcal{F}_{-\infty}^u \right)$ and $L^2\left(\mathcal{G}^{\infty}_u \right)\subset L^2\left(\mathcal{F}^{\infty}_u \right)$ . Thus,
\begin{equation*}
    \sup_{t}\sup_{f\in L^2\left(\mathcal{G}_{-\infty}^t\right), g\in L^2\left(\mathcal{G}^{\infty}_{t+k}\right) }\left|\mathrm{corr}\left(f,g\right)\right| \le \sup_{t}\sup_{f\in L^2\left(\mathcal{F}_{-\infty}^t\right), g\in L^2\left(\mathcal{F}^{\infty}_{t+k}\right) }\left|\mathrm{corr}\left(f,g\right)\right|.
\end{equation*}
Therefore, $\alpha'_n\left(k\right)\le\rho'_n\left(k\right)\le\rho_n\left(k\right)\le K_1\exp \left ( -K_2k^{\gamma } \right )$. So the condition b in Lemma 1 holds.

Since $\left \| \left ( \mathbf{T}^2 \right )^{\odot 2}  \right \|\ge \left \| \left ( \mathbf{T}^{\odot 2} \right )^2  \right \|$ and $\kappa_4\ge0$, \begin{align*}
    \min diag\left(\mathbf{\Xi}\right)&=p^{-1}\left ( 8\left \| \left ( \mathbf{T}^2 \right )^{\odot 2}  \right \| +2\left ( \kappa _4 -3\right ) \left \| \left ( \mathbf{T}^{\odot 2} \right )^2  \right \|\right )\\
    &\ge2p^{-1}\left \| \left ( \mathbf{T}^2 \right )^{\odot 2}  \right \|\\
    &=2p^{-1}\left \| \mathbf{R}^{\odot 2}  \right \|>2.
\end{align*}
So the condition c in Lemma 1 holds. By Lemma 1
\begin{equation*}
    \sup_{A\in \mathcal{A}}\left | \mathbb{P}\left ( \bm{S}_p-2\bm{1}_{\frac{n\left(n-1\right)}{2}} \in A \right ) - \mathbb{P}\left ( \bm{G}\in A \right )  \right |\lesssim \frac{B^{2/3}\log\left ( \frac{n\left ( n-1 \right ) }{2}  \right )^{\frac{1+2\gamma }{3\gamma } } }{p^{\frac{1}{9} }}+\frac{B\log\left ( \frac{n\left ( n-1 \right ) }{2}  \right )^{\frac{7}{6} } }{p^{\frac{1}{9} }},
\end{equation*}
where $\bm{G}\sim N\left(0,\mathbf{\Xi} \right)$ and $\mathcal{A}$ is a class of Borel subsets in $\mathbb{R}^{\frac{n\left(n-1\right)}{2}}$. By calculating the right formula, we can find $\mathrm{RHS}\asymp \frac{\left ( \log n \right )^{\nu} }{p^{\frac{1}{9} }} $ where $\nu=\frac{2+4\gamma }{3\gamma }\vee \frac{7}{3}$. By Condition (C2), we have $\mathrm{RHS}=o\left(1\right)$.\\
Notice that $M^2_n=\sqrt{p}\max_{1\le i<j\le n}\bm{S}_p$, Thus\\
\begin{equation*}
    \sup_{t}\left | \mathbb{P}\left ( M^2_n-2p \le t \right ) - \mathbb{P}\left ( \max_{1\le i<j\le n}\bm{G} \le t \right )  \right |=o\left(1\right).
\end{equation*}
Let $t=b_p^{1/2}x$ , $A_n=\max_{1\le i<j\le n}b_p^{1/2}\bm{G}$ , Thus
\begin{equation*}
    \sup_{x}\left | \mathbb{P}\left ( \frac{M^2_n-2p}{\sqrt{pb_p}} \le x \right ) - \mathbb{P}\left ( A_n \le x \right )  \right |=o\left(1\right).
\end{equation*}

Step 3: \\
Set $I=\left \{ \left ( i,j \right ) ;1\le i <j\le n \right \} $ for $\alpha= \left ( i,j \right ) \in I$
define $$B_\alpha=\left \{ \left ( k,l \right )\in I ;\left \{ k,l \right \}\cap \left \{ i,j \right \}\ne\emptyset ,\left ( k,l \right )\ne\alpha  \right \}.$$ Let $\eta_{\left ( i,j \right )}$ be the term corresponding to $\left ( i,j \right )$ in $b_p^{-1/2}\bm{G}$. We can find that $\eta_{\left ( i,j \right )}\sim N\left(0,1\right)$ and the covariance structure is $\mathbf{\Sigma}_p$.\\
Set $s_n=\sqrt{4\log n-\log \log n+x}$ By Lemma 2, we have\\
\begin{equation*}
\left | \mathbb{P}\left ( \max_{\alpha \in I}\eta_\alpha \le s_n \right ) -e^{-\lambda _p} \right | \le u_1+u_2,
\end{equation*}
where
\begin{align*}
    \lambda_p&=\sum_{\alpha \in I}\mathbb{P}\left (  \eta_\alpha>s_n \right ) =\frac{n\left ( n-1 \right ) }{2} \left ( 1-\Phi \left ( s_n \right )  \right ) , \\
    u_1&=\sum_{\alpha \in I}\sum_{\beta\in B_\alpha } \mathbb{P}\left (  \eta_\alpha>s_n \right )\mathbb{P}\left (  \eta_\beta>s_n \right ) =\frac{n\left ( n-1 \right ) }{2}\cdot \left ( 2\left ( n-2 \right )  \right )  \left ( 1-\Phi \left ( s_n \right )  \right ) ^2 ,\\
    u_2&=\sum_{\alpha \in I}\sum_{\beta\in B_\alpha } \mathbb{P}\left (  \eta_\alpha>s_n ,  \eta_\beta>s_n \right ) =\frac{n\left ( n-1 \right ) }{2}\cdot \left ( 2\left ( n-2 \right )  \right ) \cdot  \mathbb{P}\left (  \eta_\alpha>s_n ,  \eta_\beta>s_n \right ).
\end{align*}
If we replace $a_p$ with $b_p $ in Proof of theorem 2.1, then the proof here would be almost identical. The only difference lies in the method of verifying whether $u_2$ converges. Below, we provide a proof of the convergence of $u_2$:

Notice that $\left ( \eta_\alpha,\eta_\beta  \right )^T \sim N\left ( 0,\begin{pmatrix}
  1& \rho_p \\
 \rho_p &1
\end{pmatrix} \right ) $ and $\rho_p\le \frac{1}{3}$ we can calculate\\
\begin{align*}
    u_2&\le\frac{n\left ( n-1 \right ) }{2}\cdot \left ( 2\left ( n-2 \right )  \right ) \cdot  \mathbb{P}\left (  \eta_\alpha+  \eta_\beta>2s_n \right )\\
    &\le\frac{n\left ( n-1 \right ) }{2}\cdot \left ( 2\left ( n-2 \right )  \right ) \cdot  \mathbb{P}\left (  N\left ( 0,2+2\rho_p  \right ) >2s_n \right )\\
    &=n^3\cdot\left ( 1-\Phi \left ( \frac{\sqrt{2}s_n }{\sqrt{1+\rho_p}}  \right )  \right ) \cdot\left ( 1+o\left ( 1 \right )  \right )\\
    &\le n^3\cdot\left ( 1-\Phi \left ( \frac{\sqrt{6}s_n }{2}  \right )  \right ) \cdot\left ( 1+o\left ( 1 \right )  \right )\\
    & \asymp \frac{n^3}{s_n} e^{-\frac{3s_n^2}{4} }=O\left ( s_n^{-1} \right ) \rightarrow 0.
\end{align*}

\end{document}